\documentclass[12pt]{article}
\usepackage{amsmath}
\usepackage{latexsym}
\usepackage{amssymb}
\newtheorem{thm}{Theorem}[section]
\newtheorem{la}[thm]{Lemma}
\newtheorem{Defn}[thm]{Definition}
\newtheorem{Remark}[thm]{Remark}
\newtheorem{Note}[thm]{Note}
\newtheorem{prop}[thm]{Proposition}

\newtheorem{Example}[thm]{Example}

\newtheorem{Examples}[thm]{Examples}
\newtheorem{Problems}[thm]{Problems}

\newtheorem{problem}{Problem}
\newtheorem{Number}[thm]{\!\!}
\newenvironment{defn}{\begin{Defn}\rm}{\end{Defn}}

\newenvironment{example}{\begin{Example}\rm}{\end{Example}}

\newenvironment{rem}{\begin{Remark}\rm}{\end{Remark}}
\newenvironment{numba}{\begin{Number}\rm}{\end{Number}}
\newenvironment{proof}{{\noindent\bf Proof.}}%
                  {\nopagebreak\hspace*{\fill}$\Box$\medskip\medskip\par}   
\newcommand{\Punkt}{\nopagebreak\hspace*{\fill}$\Box$}

\newcommand{\wb}{\overline}
\newcommand{\ve}{\varepsilon}
\newcommand{\at}{\symbol{'100}}

\newcommand{\n}{\rm}

\newcommand{\mto}{\mapsto}

\newcommand{\isom}{\cong}

\newcommand{\N}{{\mathbb N}}
\newcommand{\R}{{\mathbb R}}

\newcommand{\Z}{{\mathbb Z}}

\newcommand{\F}{{\mathbb F}}
\newcommand{\bO}{{\mathbb O}}
\newcommand{\K}{{\mathbb K}}

\newcommand{\Aut}{\mbox{\n Aut}}

\newcommand{\sub}{\subseteq}

\newcommand{\GL}{\mbox{\rm GL}}

\newcommand{\id}{\mbox{\n id}}

\newcommand{\dsemi}{\mbox{$\times\!$\rule{.15mm}{2mm}}\,}
\begin{document}
\begin{center}
{{\Large\bf Smooth Lie groups over local fields of\vspace{2mm}
positive characteristic need not be analytic}}\\[6mm]
{\bf Helge Gl\"{o}ckner}\vspace{5mm}
\end{center}
\begin{abstract}
\hspace*{-7.2 mm}
We describe finite-dimensional smooth Lie groups over local fields of positive characteristic
which do not admit an analytic Lie group structure
compatible with the given topological group structure,
and $C^n$-Lie groups without a compatible
$C^{n+1}$-Lie group structure,
for each $n\geq 1$.
We also present examples of non-analytic, smooth automorphisms,
as well as $C^n$-automorphisms which fail to be $C^{n+1}$.\vspace{2mm}
\end{abstract}
{\footnotesize Subject classification:
22E20 (main); 
20G25, 
22D05, 
22D45, 
26E30.\\[2mm] 
}\\[5mm]
\begin{center}
{\large\bf Introduction}
\end {center}
It is a classical fact (the core of which goes back to
Schur, \cite{Shu}) that every finite-dimensional
real Lie group of class $C^k$
(where $k\in \N\cup\{\infty\}$)
admits a $C^k$-compatible analytic Lie group structure;
see
\cite[\S4.4]{MaZ}
and the references given there for further information.
Based on Lazard's characterization of
analytic $p$-adic Lie groups within the class
of locally compact groups \cite{Laz},
it was shown in \cite{ANA}
that also every finite-dimensional
$C^k$-Lie group over a local field
of characteristic $0$
admits a $C^k$-compatible
$\K$-analytic Lie group structure.
Here, $C^k$-maps and Lie groups
are understood in the sense of~\cite{DIF},
where a setting of differential
calculus over arbitrary topological fields
was developed. A mapping between finite-dimensional vector spaces
over a local field is $C^1$ if and only if it is strictly differentiable at
each point~\cite[Rem.\,5.4]{IMP},
where strict differentiability is a classical concept, \cite[1.2.2]{Bo1}.\\[3mm]
In the present article, we show by example
that the situation is entirely different
in the case of positive characteristic.
We prove:\\[4mm]
{\bf Main Theorem.}
{\em Let $\,\K$ be a local field of positive characteristic.
Then\/}:
\begin{itemize}
\item[(a)]
{\em There exists a one-dimensional smooth $\K$-Lie group
which does not admit a $\K$-analytic
Lie group structure compatible with its topological group
structure.}
\item[(b)]
{\em For each positive integer $n$, there exists
a one-dimensional $\K$-Lie group of class~$C^n$
which does not admit a $\K$-Lie group structure of class $C^{n+1}$
compatible with the given topological group structure.}\vspace{2mm}
\end{itemize}
(See Theorem~\ref{mainthm}
and Theorem~\ref{mainthm2}).
We also describe a $C^n$-automorphism of a finite-dimensional
$\K$-analytic Lie group
which is not $C^{n+1}$ (Proposition~\ref{mainwork1}),
and a non-analytic $C^\infty$-automorphism (Proposition~\ref{mainwork2}).\\[3.1mm]
As a tool, we use concepts from George Willis' structure theory
of totally disconnected, locally compact groups as developed
in \cite{Wi1}, \cite{Wi2} and analyzed for Lie groups over local fields
in \cite{SCA} and \cite{POS}.\\[3.1mm]
While it is quite justified to restrict attention
to analytic Lie groups in the real
and $p$-adic cases, our results make it clear
that this habit is inappropriate in the
case of positive characteristic.
Because genuinely larger classes of groups
can be captured in this way,
it may very well be of interest
to study Lie groups
over such fields which are merely smooth,
or merely~$C^n$. It is not hopeless to prove substantial results
in the absence of analyticity.
For example, it is possible to discuss
contraction groups (opposite subgroups), tidy subgroups
and the scale
function for
automorphisms of mere $C^1$-Lie groups, \cite{POS}.
\section{Preliminaries and Notation}
%
%
We recall some definitions and basic
facts from differential calculus (and
linear algebra) over local fields,
and set up our notational conventions.
\begin{numba}
Let $E$ and $F$ be (Hausdorff) topological
vector spaces
over a non-discrete topological field~$\K$,
$U\sub E$ be open, and $f\!: U\to F$ be a map.
Following \cite{DIF},
we call $f$ a $C^1$-map
if it is continuous and there
exists a (necessarily unique)
continuous map $f^{[1]}\!: U^{[1]}\to F$
on $U^{[1]}:=\{(x,y,t)\in U\times E\times \K\!:
x+ty\in U\}$ such that $f^{[1]}(x,y,t)=\frac{1}{t}(f(x+ty)-f(x))$
for all $(x,y,t)\in U^{[1]}$ with $t\not=0$.
Then $f'(x)\!: E\to F$, $y\mto f'(x).y:=f^{[1]}(x,y,0)$ is a continuous linear
map, for each $x\in U$ \cite[Prop.\,2.2]{DIF}.
Inductively, $f$ is called $C^{n+1}$ if it is
$C^1$ and $f^{[1]}$ is $C^n$;
it is  $C^\infty$ (or {\em smooth}) if it is $C^n$
for all~$n$.
\end{numba}
\begin{numba}
For
maps between open subsets
of finite-dimensional
real vector spaces,
the present definitions
of $C^n$-maps are equivalent to the usual ones
(cf.\ \cite[Prop.\ 7.4 and 7.7]{DIF}).
\end{numba}
\begin{numba}\label{deflie}
Compositions of composable $C^r$-maps
are $C^r$ \cite[Prop.\,4.5]{DIF},
the Chain Rule holds \cite[Prop.\,3.1]{DIF},
and being $C^r$ is a local property,
for each $r\in \N\cup\{\infty\}$
\cite[La.\,4.9]{DIF}.
This makes it possible to define
$C^r$-manifolds modeled on a topological
$\K$-vector space~$E$,
in the usual way.
A {\em $\K$-Lie group of class $C^r$}
is a group, equipped with a $C^r$-manifold structure
modeled on a topological $\K$-vector space, which
makes the group operations~$C^r$.
If $G$ is a $C^r$-Lie group, we
write $L(G):=T_1(G)$ for its (geometric) tangent space at the identity.\footnote{If
$r\geq 3$, then $L(G)$ has a natural Lie algebra structure,
but we shall not use it.}
Given a homomorphism $\alpha\!: G\to H$ which is at least~$C^1$,
as usual we abbreviate $L(\alpha):=T_1(\alpha)$.
\end{numba}
\begin{numba}
For the definition of
$\K$-analytic maps between
open subsets of Banach spaces over a complete
valued field~$\K$, we refer to~\cite{Bo1} (see also
\cite{Ser} for the finite-dimensional case).
Any $\K$-analytic map is $C^\infty$ by \cite[Prop.\,7.20]{DIF}.
Hence every $\K$-analytic Lie group modeled
on a Banach space can also be considered as
a $C^\infty$-Lie group in the sense of {\bf \ref{deflie}}.
\end{numba}
For functions of a single variable,
an equivalent description of the $C^n$-property
is available, which is much easier to work with.
This is essential for
our constructions of pathological homomorphisms and Lie groups,
which necessitate to verify by hand that certain mappings
defined on open subsets of the ground field are~$C^n$.
It is also clear from the equivalent description that a map
between open subsets of an ultrametric field $\K$ is $C^n$ in our sense if and only if it
is so in the sense of~\cite{Sch}.
\begin{la}\label{easierdiff}
{\rm(cf.\ \cite[Prop.\,6.9]{DIF}).}
Let $U\sub \K$ be an open subset of a topological field~$\K$
and $F$ be a topological $\K$-vector space.
Given a map $f\!: U\to F$, we set
$f^{<0>}:=f$. Then the following holds:
\begin{itemize}
\item[\rm (a)]
The map $f$ is $C^1$
if and only if it is continuous and there exists
a $($necessarily unique$)$
continuous map $f^{<1>}\!: U\times U\to F$
such that
\[
f^{<1>}(x_1,x_2)\; =\; \frac{1}{x_1-x_2}(f(x_1)-f(x_2))
\]
for all $x_1,x_2\in U$ such that $x_1\not=x_2$.
\item[\rm (b)]
Having characterized $C^k$-maps
and defined $f^{<k>}\!: U^{k+1}\to F$
for all $k=0,\ldots, n-1$, where $n\in \N$,
we have:
$f$ is~$C^n$ if and only if $f$ is
$C^{n-1}$
and there exists a continuous map
$f^{<n>}\!: U^{n+1}\to F$
such that
\begin{eqnarray*}
\lefteqn{f^{<n>}(x_1,x_2,\ldots, x_{n+1})}\\
\!\!\!\!\!&\!\!\!=&\!\!\!\!\frac{1}{x_1-x_2}\left(f^{<n-1>}(x_1,x_3,\ldots, x_{n+1})
-f^{<n-1>}(x_2,x_3,\ldots, x_{n+1})\right)
\end{eqnarray*}
for all $x_1,\ldots,x_{n+1}\in U$
such that $x_1\not=x_2$.
\end{itemize}
Here $f^{<n>}$ is uniquely determined
and is symmetric in its $n+1$ arguments.
Furthermore, $n!\, D^nf(x)=\frac{d^nf}{dx^n}(x)=:f^{(n)}(x)$
for all $x\in U$, where $D^nf(x):=f^{<n>}(x,\ldots, x)$.\Punkt
\end{la}
For further
information on differential calculus over topological fields
and the corresponding manifolds
and Lie groups,
we refer to \cite{DIF} 
and \cite{ANA}--\cite{ZOO}.
\begin{numba}
Recall that a {\em local field\/}
is a totally disconnected, locally compact, non-discrete topological field.
On any local field, there exists a {\em natural absolute value\/}
$|.|\!: \K\to [0,\infty[$,
defined via $|0|:=0$ and $|x|:=\Delta_\K(\lambda_x)$
for $x\not=0$, where $\lambda_x\!: \K\ \to \K$, $\lambda_x(y):=xy$
is multiplication by~$x$ and $\Delta_\K\!: \Aut(\K,+)\to\; ]0,\infty[$
the modular function.
Then $|\K^\times|=\langle q\rangle\leq (]0,\infty[,\cdot)$
for some prime power~$q$,
called the {\em module of~$\K$}.
For further information on local fields,
see~\cite{Wei}.
\end{numba}
The following lemma provides information concerning a special type
of linear maps which will play a crucial role in our constructions.
\begin{la}\label{crucialla}
Let $V$ be a finite-dimensional
vector space over a local field $\K$
with module~$q$.
Let $|.|\!: \wb{\K}\to[0,\infty[$
be the unique extension of the natural absolute value
on~$\K$ to an absolute value on the algebraic closure $\wb{\K}$ of~$\K$.
Let $\alpha\in\GL(V)$ be a $\K$-linear automorphism
of~$V$, with eigenvalues $\lambda_1,\ldots, \lambda_d$ in~$\wb{\K}$, where $d:=\dim_\K(V)$.
Suppose that $|\lambda_j|>1$ for each $j\in \{1,\ldots, d\}$
and $\prod_{j=1}^d |\lambda_j|=q$.
Consider $V$ as a $\K[T]$-module
via $\big(\sum_{j=0}^na_jT^j\big).v:=\sum_{j=0}^na_jA^jv$.
Then
\begin{itemize}
\item[\rm (a)]
$V$ is an irreducible $\K[T]$-module;
\item[\rm (b)]
$|\lambda_j|=q^{1/d}$ for each $j\in \{1,\ldots, d\}$; and
\item[\rm (c)]
There exists an ultrametric norm $\|.\|$ on~$V$
such that $\|\alpha(x)\|=q^{1/d}\|x\|$
for each $x\in V$.
\end{itemize}
\end{la}
\begin{proof}
(a) Let $\{0\}\not=U\sub V$ be a $\K[T]$-submodule
of $V$; we show that $U=V$.
To this end, let $e:=\dim_\K(U)\in \{1,\ldots, d\}$
and note that $\lambda_{j_1},\ldots, \lambda_{j_e}$ are
the eigenvalues of $\alpha|_U^U$ in $\wb{\K}$,
for suitable integers $1\leq j_1<j_2<\cdots<j_e\leq d$.
Set $J:=\{j_k\!: k=1,\ldots, e\}$ and $r:=\prod_{j\not\in J}|\lambda_j|$.
Then $r\geq 1$, with equality if and only if $J=\{1,\ldots,d\}$
(and hence $U=V$).
But
\[
q\,=\, |\det(\alpha)|\,=\, \prod_{j=1}^d |\lambda_j|\,=\, r\, \prod_{j\in J}|\lambda_j|\,,
\]
where $1<\prod_{j\in J}|\lambda_j|=|\det(\alpha|_U^U)|\in q^\Z$
and thus $\prod_{j\in J}|\lambda_j|\geq q$.
Hence  $\prod_{j\in J}|\lambda_j|=q$ and $r=1$,
entailing that $U=V$.\vspace{2mm}

(b) For each $j$, let $(V_{\wb{\K}})_{\lambda_j}:=\{ v\in V_{\wb{\K}}\!: (\alpha_{\wb{\K}}-\lambda_j)^dv=0\}$
be the generalized eigenspace of $\alpha_{\wb{\K}}$
in $V_{\wb{\K}}$ corresponding to~$\lambda_j$,
where $V_{\wb{\K}}$ is the $\wb{\K}$-vector space obtained from~$V$ by extension of scalars.
Let $R:=\{|\lambda_j|\!: j\in\{1,\ldots, d\}\}$;
given $\rho\in R$ define $J_\rho:=\{j\in \{1,\ldots, d\}\!: |\lambda_j|=\rho\}$
and $W_\rho:=\bigoplus_{j\in J_\rho} (V_{\wb{\K}})_{\lambda_j}$.
By \cite[Ch.\,II, (1.0), p.\,81]{Mar}, each $W_\rho$ is defined over~$\K$, i.e.,
$W_\rho= (V_\rho)_{\wb{\K}}$.
Note that $V_\rho:=W_\rho\cap V$
is $\alpha$-invariant and hence a $\K[T]$-submodule of~$V$.
Thus
$V=\bigoplus_{\rho\in R} V_\rho$
is a direct sum of $\K[T]$-submodules.
By (a), such a decomposition is only possible if $R=\{\rho\}$ is a singleton
and thus $V=V_\rho$. Now $q=|\lambda_1|\cdots|\lambda_d|=\rho^d$
entails $\rho=q^{1/d}$.\vspace{2mm}

(c) We choose a $\wb{\K}$-basis $v_1,\ldots, v_d$ of $V_{\wb{\K}}$
with respect to which $\alpha_{\wb{\K}}$ has Jordan canonical
form, and let $\|.\|$ be the maximum norm
on $V_{\wb{\K}}$ with respect to this basis.
Since $q^{1/d}>1$ and $|1|=1$, the ultrametric inequality entails that
$\|\alpha_{\wb{\K}}(v_j)\|=q^{1/d}\|v_j\|$
for each $j\in\{1,\ldots, d\}$ (cf.\
\cite[p.\,209, Eqn.\,(1)]{SCA}).
Hence $\|\alpha_{\wb{\K}}(x)\|=q^{1/d}\|x\|$ for each $x\in V_{\wb{\K}}$,
using that $\|.\|$ is a maximum norm.
Now the restriction of $\|.\|$
to $V$ is norm with the desired property.\vspace{-2mm}
\end{proof}
\begin{numba}
{\bf Conventions.}
Throughout the article, $q$ is a prime power,
$\F_q$ a finite field
of cardinality~$q$ and $\K:=\F_q(\!(X)\!)$
the field of formal Laurent series over~$\F_q$,
unless said otherwise.\footnote{Up to isomorphism, every
local field of positive characteristic is of this form~\cite{Wei}.}
We let $\bO:=\F_q[\![X]\!]$ be the valuation ring of~$\K$
(the ring of formal power series over~$\F_q$).
We let
$|.|\!: \K\to [0,\infty[$ be the natural
absolute value on $\K$, given by $|\sum_{k=n}^\infty a_k X^k|:=q^{-n}$
for $a_n,a_{n+1},\ldots \in \F_q$ with $a_n\not=0$.
We fix an algebraic closure $\wb{\K}$ and use the same symbol, $|.|$,
for the unique extension of the absolute value $|.|$ on $\K$ to
an absolute value on~$\wb{\K}$.
We use $\R$ and $\Z$ in their usual meanings and write $\N:=\{1,2,\ldots\}$,~$\N_0:=\N\cup\{0\}$.\\[3mm]
In the following, by a ``Lie group''
we always mean a Lie group
modeled on a finite-dimensional vector space.
A {\em $C^r$-automorphism\/} of a $C^r$-Lie group~$G$
is an automorphism $\alpha\!: G\to G$
such that both $\alpha$ {\em and its inverse\/}
are $C^r$.
The same convention applies to
$\K$-analytic automorphisms,
$\K$-analytic diffeomorphisms and $C^r$-diffeomorphisms.
Given a topological group~$G$,
its group of bicontinuous automorphisms will
be denoted by $\Aut(G)$.
If $(E,\|.\|)$ is a normed
$\K$-vector space, we set
$B_r^E:=\{x\in E\!: \|x\|<r\}$ for $r>0$;
if $E$ is understood, we abbreviate $B_r:=B_r^E$.
\end{numba}
\section{$\!\!\!$Criteria making a map {\boldmath $C^n$}, or preventing it}
In this section, we describe simple criteria
ensuring that a mapping is $C^n$,
resp., that it is not $C^n$.
\begin{la}\label{tandif}
Let $n\in \N$ and $f\!: U \to\K$ be a map,
defined on an open subset $U\sub \K$.
Suppose that there exist real numbers $\theta >n$
and $b > 0 $ such that
\begin{equation}\label{highvan}
|f(y)-f(x)|\leq b \, |y-x|^\theta \qquad \mbox{for all $\, x,y \in U$.}
\end{equation}
Then $f$ is $C^n$
and $D^kf(x)=0$
for all $k\in \{1,\ldots,n\}$ and $x\in U$.
\end{la}
\begin{proof}
This is a special case of \cite[Thm.\,29.12]{Sch}.
\end{proof}
\begin{la}\label{nonCn}
Let $E\not=\{0\}$ and $F$ be
normed $\K$-vector spaces
and $f\!: U\to F$ be a mapping on an open subset $U\sub E$.
Let $n\in \N_0$.
Suppose that there exist $x\in U$, $\theta\in \; ]n,n+1[$ and positive real numbers
$a,b$ such that
\begin{equation}\label{badasympt}
a\, \|y-x \|^\theta \;\leq \; \|f(y)-f(x)\|\;\leq\;
b\, \|y-x\|^\theta
\end{equation}
for all $y\in U$.
Then $f$ is not of class $C^{n+1}$.
\end{la}
\begin{proof}
There exists $0\not=y\in E$ such that $x+\bO y\sub U$.
Then $g\!: \bO\to F$, $g(t):=f(x+ty)-f(x)$
satisfies
\begin{equation}\label{badas2}
A\, |t|^\theta \leq \|g(t)\|\leq B\, |t|^\theta\qquad \mbox{for all $t\in \bO$,}
\end{equation}
with $A:= a\, \|y\|^\theta$, $B:=b\, \|y\|^\theta$. It suffices to show that~$g$ is not $C^{n+1}$.
To reach a contradiction, assume that $g$ is $C^{n+1}$.
By \cite[Thm.\,5.1]{DIF}, $g$ admits a
Taylor expansion of order $n+1$ around $0$:
\[
g(t)\; =\; \sum_{j=0}^{n+1} a_j t^j\, + \, t^{n+1}\rho(t)\qquad \mbox{for all $t\in \bO$,}
\]
where $\rho\!:\bO\to F$ is continuous and satisfies $\rho(0)=0$.
But such an expansion is incompatible with the asymptotic described in~(\ref{badas2}).
To see this, note first that $a_k=0$ for $k=0,\ldots, n$.
Indeed,
$g(0)=0$ entails that $a_0=0$, and if we already know that $a_0=\cdots =a_{k-1}=0$,
letting $t\to 0$ in
$t^{-k}g(t)=a_k+\sum_{j=k+1}^{n+1}t^{j-k}a_j+ t^{n+1-k}\rho(t)$
we find that also $a_k=0$
(since $t^{-k}g(t)\to 0$, by (\ref{badas2})).
Hence $g(t)=(a_{n+1} +\rho(t))\, t^{n+1}$
and thus
\[
A\;\leq \; \frac{\|g(t)\|}{|t|^\theta} \; =\; \|a_{n+1}+\rho(t)\|\cdot |t|^{n+1-\theta}\, ,
\]
where the right hand side tends to $0$ as $t\to 0$. This contradicts $A>0$.
\end{proof}
\section{Automorphisms with specific properties}\label{secauto}
In this section, we discuss a family of automorphisms of the additive
topological group $\bO:=\F_q[\![X]\!]$, some of which have specific
pathological properties.
\begin{numba}\label{defalpbet}
Given a strictly ascending sequence $\lambda:=(\ell_j)_{j\in \N_0}$
of positive integers $\ell_0<\ell_1<\ell_2< \cdots$,
we define a mapping
$\beta\!: \bO\to \bO$ via
\[
\beta\Big(
\sum_{j=0}^\infty a_j X^j\Big)\; := \;
\sum_{j=0}^\infty a_j X^{\ell_j}\qquad
\mbox{for $\;(a_j)_{j\in \N_0}\in \F_q^{\N_0}$.}
\]
We define $\alpha \! : \bO\to \bO$ via $\alpha(z):= z+ \beta(z)$.
\end{numba}
\begin{la}\label{basalpbet}
The map $\beta$ from {\bf \ref{defalpbet}} is a continuous homomorphism,
and $\alpha$ is an automorphism of the topological group $(\bO, +)$.
\end{la}
\begin{proof}
Clearly $\beta$ and $\alpha$ are homomorphisms
of groups. Furthermore, $\beta$ is continuous,
since $\bO=\F_q^{\N_0}$ is equipped with the product
topology and each coordinate of $\beta(z)$
is either constant or only depends on one coordinate of~$z$.\\[3mm]
To see that $\alpha$ is a bijection,
let $z=\sum_{j=0}^\infty a_jX^j$ and $w\in \sum_{j=0}^\infty b_j X^j$
be elements of $\bO$, where $a_j,b_j\in \F_q$.
Equating coefficients, we see that $\alpha(z)=w$
holds if and only if $a_0=b_0$ and
\begin{equation}\label{splitcase}
\!\!\!\!\!(\forall i\in \N)\qquad a_i\; =\;
\left\{
\begin{array}{cl}
b_i & \; \mbox{if $\;i\not=\ell_j\,$ for each $j \in \N_0$;}\\
b_i-a_j & \; \mbox{if $\;i=\ell_j\,$ for some $j\in \N_0$.}
\end{array}
\right.
\end{equation}
Since $j<i$ in the second case of (\ref{splitcase}), given $w$
we can solve for $a_0$, $a_1$, $\ldots$ in turn
and deduce that $z\in \bO$ with $\alpha(z)=w$ exists
and is uniquely determined. Thus $\alpha$ is an isomorphism of groups.
It is clear from the preceding formula that, for fixed $i$,
the $i$th coordinate $a_i$ of $z=\alpha^{-1}(w)$ only depends
on finitely many coordinates,
$b_0, b_1,\ldots, b_i$, of $w$. Hence $\alpha^{-1}$ is continuous.
\end{proof}
\begin{la}\label{lamgivesas}
If there exist real numbers
$\theta>0$ and $K>0$ such that
\begin{equation}\label{condlam}
j\theta- K\; \leq \; \ell_j\; \leq \; j \theta+ K\qquad\mbox{for all $\,j\in \N_0$,}
\end{equation}
then there exist positive real numbers $a<b$ such that
\[
a\, |y-x|^\theta \; \leq \; |\beta(y)-\beta(x)|
\;\leq \; b\, |y-x|^\theta\, ,\quad \mbox{for all $\,x,y\in \bO$.}
\]
\end{la}
\begin{proof}
We claim that the assertion holds with $a:=\min\{q^{-\ell_0}, q^{-K}\}$ and $b:=q^K$.
Since $\beta(y)-\beta(x)=\beta(y-x)$, we may assume without loss of generality that $x=0$.
We distinguish three cases.\\[3mm]
1. If $y=0$, then the assertion is trivial.\\[3mm]
2. If $|y|=1$, then $|\beta(y)|=q^{-\ell_0}\geq a=a|y|^\theta$
and $b|y|^\theta=b\geq 1\geq q^{-\ell_0}=|\beta(y)|$.\\[3mm]
3. If $|y|=q^{-j}$ with $j\in \N$, then $\ell_j/j=\theta +\ve_j$ for some
$\ve_j\in \big[{-K/j}, K/j\big]$, by~(\ref{condlam}).
Thus
\[
|\beta(y)|=q^{-\ell_j}=(q^{-j})^{\ell_j/j}=|y|^{\ell_j/j}=|y|^{\theta+\ve_j}=|y|^\theta |y|^{\ve_j}\, ,
\]
where $|y|^{\ve_j}=q^{-j\ve_j}$ with $j\ve_j\in [{-K},K]$
and so $|y|^{\ve_j}\in [q^{-K},q^K]$.
Hence indeed $a|y|^\theta\leq q^{-K}|y|^\theta\leq|\beta(y)|\leq q^K|y|^\theta=b|y|^\theta$.
\end{proof}
The following proposition provides $C^n$-automorphisms of the $1$-dimensional
$\K$-analytic Lie group $(\bO,+)$ which are not $C^{n+1}$.
\begin{prop}\label{mainwork1}
Given $n\in \N$, choose $\theta\in \;]n, n+1[$
and define $\ell_j:=1+[ j\theta]$ for each $j\in \N_0$,
where $[.]$ is the Gau\ss{} bracket $($integer part$)$.
Let $\alpha\!: \bO\to \bO$ be as in {\bf\ref{defalpbet}}.
Then $\alpha$ is a $C^n$-automorphism of $(\bO,+)$
which is not a $C^{n+1}$-map.
\end{prop}
\begin{proof}
Since $j\theta\leq \ell_j\leq j\theta+1$ for each $j\in \N_0$, Lemma~\ref{lamgivesas}
provides real numbers $0<a<b$ such that $a\,|y-x|^\theta\leq |\beta(y)-\beta(x)|\leq b\, |y-x|^\theta$
for all $x,y\in \bO$. Hence $\beta$ is $C^n$
(by Lemma~\ref{tandif}) but not $C^{n+1}$
(by Lemma~\ref{nonCn}).
Furthermore, $\beta'(x)=D^1\beta(x)=0$ for all $x\in \bO$,
by Lemma~\ref{tandif}.
Hence $\alpha=\id_\bO+\beta$ is $C^n$-map which is not $C^{n+1}$,
and such that $\alpha'(x)=1$ for each $x\in \bO$.
As $\alpha$ is an automorphism of groups by Lemma~\ref{basalpbet}
and hence a bijection, the Inverse Function Theorem
for $C^n$-maps (\cite[Thm.\,7.3 \& Rem.\,5.4]{IMP} or \cite[Thms.\,27.5 \& 77.2]{Sch})
implies that $\alpha$ is a $C^n$-diffeomorphism.
\end{proof}
To identify suitable smooth mappings as non-analytic,
we shall use the following simple fact
(see \cite[4.2.3 and 3.2.5]{Bo1}):
\begin{la}\label{givnonana}
Let $E$ and $F$ be ultrametric Banach spaces over $\K$,
$U\sub E$ be an open $0$-neighbourhood and
$f\!: U\to F$ be a $\K$-analytic function
such that
\[
\lim_{x\to 0} \|f(x)\|\, /\, \|x\|^n\; =\; 0\qquad\mbox{for each $\, n\in \N_0$.}
\]
Then there exists a $0$-neighbourhood $W\sub U$ such that $f|_W=0$.\Punkt
\end{la}
We find it convenient to formulate another preparatory lemma.
\begin{la}\label{preparat}
In {\bf \ref{defalpbet}}, choose $\lambda=(\ell_j)_{j\in \N_0}$
such that $\, \ell_j/j\to \infty$ as $\,j\to\infty$.
Then, for each $n\in \N$, there exists a real number $c_n>0$
such that
\begin{equation}\label{givinftan}
|\beta(y)-\beta(x)|\;\leq\; c_n\, |y-x|^n\qquad\mbox{for all $\,x,y\in \bO$.}
\end{equation}
\end{la}
\begin{proof}
Let $n\in \N$. Since $\ell_j/j\to\infty$,
there exists $k\in\N$ such that $\ell_j/j\geq n$ for all
$j\geq k$. We claim that (\ref{givinftan})
holds with $c_n:=q^{kn}$.
It suffices to verify
(\ref{givinftan}) for $x=0$, and we may assume that $y\not=0$
(cf.\ proof of Lemma~\ref{lamgivesas}).
If $|y|=q^{-j}$ with $j\geq k$,
then $|\beta(y)|=q^{-\ell_j}=(q^{-j})^{\ell_j/j}\leq (q^{-j})^n
=|y|^n\leq c_n\,|y|^n$,
as required.
If $|y|=q^{-j}$ with $j\in\{0,1,\ldots, k-1\}$,
then $|\beta(y)|\leq 1 \leq q^{kn}q^{-jn}=c_n\, |y|^n$
as well.
\end{proof}
We now provide examples of smooth, non-analytic automorphisms
of the $1$-dimensional $\K$-analytic Lie group~$(\bO,+)$.
\begin{prop}\label{mainwork2}
In {\bf \ref{defalpbet}}, choose $\lambda=(\ell_j)_{j\in \N_0}$
such that $\, \ell_j/j\to \infty$ as\linebreak
$j\to\infty$.
Then the following holds:
\begin{itemize}
\item[\rm (a)]
$\alpha$ is a $C^\infty$-automorphism of $(\bO,+)$.
\item[\rm (b)]
$\alpha$ is not $\K$-analytic.
%
%
\end{itemize}
\end{prop}
\begin{proof}
(a) Lemmas~\ref{preparat} and \ref{tandif} show that
$\beta$ is smooth, with
\begin{equation}\label{prpbet}
D^j\beta(x)\,=0\,\qquad\mbox{for all $\,j\in\N$ and $\,x\in\bO$.}
\end{equation}
As in the proof of Proposition~\ref{mainwork1},
we infer that $\alpha$ is a $C^\infty$-automorphism.\vspace{2mm}

(b) Since $\beta=\alpha-\id_\bO$ where $\id_\bO$ is $\K$-analytic,
it suffices to show that $\beta$ is not $\K$-analytic.
However, if $\beta$ was $\K$-analytic, then $\beta$
would vanish on some $0$-neighbourhood,
by (\ref{givinftan}) and Lemma~\ref{givnonana}.
As $\beta$ is injective, we have reached a contradiction.
\end{proof}
Note that functions $f$ of one variable
with $D^jf=0$ for small (or all) positive integers
are familiar objects in non-archimedian
analysis.
Examples of such functions in the $p$-adic case
(which are {\em not\/} homomorphisms)
are described (or suggested) in Example~26.4 as well as Exercises 28.F, 29.F and 29.G
in~\cite{Sch}.
Our definition of $\beta$ in {\bf\ref{defalpbet}}
was stimulated by these examples.
%
%
%
%
%
%
\section{Calibrations obtained from contractions}
In this section,
we provide techniques which shall enable
us
to check that certain Lie groups
constructed with the help of the
pathological automorphisms from Section~\ref{secauto}
do not admit better-behaved Lie group structures
compatible with their given topological group structure.
The basic idea is that suitable contractive automorphisms
of $\K$-Lie groups determine a ``calibration''
which can be used to estimate the norms
of elements close to the identity
(with respect to a given chart taking $1$ to $0$).
As a tool, we shall employ concepts
from the structure theory
of totally disconnected, locally compact groups
developed in \cite{Wi1}, \cite{Wi2}
and analyzed for
Lie groups over local fields in~\cite{SCA} and~\cite{POS}.
We briefly recall the relevant definitions.
\begin{numba}
Let $G$ be a totally disconnected, locally compact topological group
and $\alpha$ be a bicontinuous automorphism of~$G$.
Then there exists a compact open subgroup $U\sub G$ which is {\em tidy for $\alpha$},
meaning that the following holds:
\begin{itemize}
\item[{\bf T1}]
$U=U_+U_-$, where $U_\pm:=\bigcap_{k\in \N_0} \;\alpha^{\pm k}(U)\,$; and
\item[{\bf T2}]
The subgroup $U_{++}:=\bigcup_{k\in \N_0}\; \alpha^k(U_+)\,$ is closed in $G$.
\end{itemize}
%
Here $\alpha(U_+)\supseteq U_+$ and $\alpha(U_-)\sub U_-$.
The index $s_G(\alpha):=[\alpha(U_+)\,:\, U_+]$ is finite and is independent
of the tidy subgroup; it is called the {\em scale of $\alpha$}.
If $U$ is tidy for~$\alpha$, then it is also tidy for $\alpha^{-1}$.
(See \cite{Wi1} and \cite{Wi2}).
\end{numba}
Let us say that a bicontinuous automorphism $\alpha\!: G\to G$ is
{\em contractive\/} if $\alpha^n(x)\to 1$ as $n\to\infty$, for each $x\in G$.
A sequence $(U_n)_{n\in \N_0}$
of identity neighbourhoods is called a {\em fundamental sequence\/}
if $U_0\supseteq U_1\supseteq U_2\supseteq \cdots$
and every identity neighbourhood of $G$ contains some~$U_n$.
\begin{numba}\label{Wang}
If $\alpha$ is contractive and $U\sub G$ a compact identity neighbourhood,
then $\{\alpha^n(U)\!: n\in \N_0\}$
is a basis for the filter of identity neighbourhoods of~$G$,
as a consequence of \cite[Prop.\,2.1]{Wan}.
\end{numba}
\begin{la}\label{wilcontra}
Let $\alpha$ be a contractive automorphism of a totally disconnected, locally compact
topological group~$G$. Then $s_G(\alpha)=1$, and the following holds:
\begin{itemize}
\item[\rm (a)]
There exists a compact, open subgroup $U\sub G$ such that $\alpha(U)\sub U$.
\item[\rm (b)]
A compact, open subgroup $U\sub G$ is tidy for $\alpha$ if and
only if $\alpha(U)\sub U$.
In this case, $U_-=U$,
\begin{equation}\label{useful}
U_+\,:=\, \bigcap_{k\in \N_0}\alpha^k(U)\,=\, \{1\}\, ,  \quad \mbox{and}\quad
U_{--}\,:=\, \bigcup_{k\in \N_0} \alpha^{-k}(U_-)\,=\,G\,.
\end{equation}
\item[\rm (c)]
$(\alpha^n(U))_{n\in \N_0}$ is a fundamental sequence of identity neighbourhoods for
$G$, for
each subgroup $U\sub G$ which is tidy for $\alpha$.
Each of the groups $\alpha^n(U)$ is tidy for~$\alpha$. Hence
$G$ has arbitrarily small subgroups tidy for~$\alpha$.
\item[\rm (d)]
If $G$ is a $C^1$-Lie group over a local field~$\K$ and $\alpha$ is a contractive $C^1$-automorphism,
let $\beta:=L(\alpha)$ be the associated $\K$-linear automorphism
of~$L(G)$.
Then $|\lambda|<1$
for every eigenvalue $\lambda$ of $\beta$ in~$\wb{\K}$.
\end{itemize}
\end{la}
\begin{proof}
(a) As tidy subgroups always exist, (a) follows from (b).\vspace{2mm}
%

(b) If $U\sub G$ is tidy for~$\alpha$,
then $\alpha^n(x)\to 1$ as $n\to\infty$
for $x\in U$ entails that $x\in U_-$ (cf.\ \cite[La.\,9]{Wi1}).
Hence $U=U_-$ and hence $\alpha(U)\sub U$.
If, conversely, $U\sub G$ is a compact open subgroup such that
$\alpha(U)\sub U$, then $\alpha^{-1}(U)\supseteq U$ and thus
$U_-=U$. Since $G$ is Hausdorff, {\bf \ref{Wang}}
shows that $U_+=\bigcap_{n\geq 0}\alpha^n(U)=\{1\}$
and hence $U_{++}=\{1\}$, which is closed.
Thus $U$ satisfies {\bf T1} and {\bf T2}.
Since $U_+=\{1\}$, furthermore $s_G(\alpha)=[\{1\}\,:\,\{1\}]=1$.
To complete the proof of (\ref{useful}),
let $x\in G$. Then $\alpha^n(x)\in U$ for some $n\in \N_0$
and thus $x\in \alpha^{-n}(U)=\alpha^{-n}(U_-)\sub U_{--}$.
Hence $U_{--}=G$.\vspace{2mm}

(c) As $\alpha^{n+1}(U)\sub \alpha^n(U)$ for each $n$,
the claims follow from (b) and {\bf \ref{Wang}}.\vspace{2mm}

(d)
As $G$ has small subgroups tidy for~$\alpha$,
results from~\cite{POS} show that
$G_\alpha:=\{x\in G : \, \mbox{$\alpha^n(x)\to 1$ as $n\to\infty$}\}$,
$G_{\alpha^{-1}}:=\{x\in G : \, \mbox{$\alpha^{-n}(x)\to 1$ as $n\to\infty$}\}$
and $G_0:=\{x\in G : \, \mbox{$\alpha^\Z(x)$ is relatively compact}\}$
are $C^1$-Lie subgroups of $G$,
and that the multiplication map
\[
G_\alpha\times G_0\times G_{\alpha^{-1}}\; \to\; G\,, \qquad
(x,y,z)\mto xyz
\]
is a $C^1$-diffeomorphism onto an open subset of~$G$.
As $G_\alpha=G$, this implies that $G_0=G_{\alpha^{-1}}=\{1\}$,
whence $L(G_0)=L(G_{\alpha^{-1}})=\{0\}$.
However, it is also shown in \cite{POS}
that $L(G_{\alpha^{-1}})_{\wb{\K}}$
(resp., $L(G_0)_{\wb{\K}}$)
is the direct sum of all the generalized eigenspaces of $\beta_{\wb{\K}}$
to eigenvalues of absolute value $>1$ (resp., of absolute value~$1$).
Hence such eigenvalues cannot exist.
\end{proof}
In view of Lemma~\ref{wilcontra}\,(b) and\,(c), the following definition makes sense:
\begin{defn}\label{defncalib}
Let $G$ be a totally disconnected, locally compact group,
$\alpha$ be a contractive automorphism of~$G$ and $U\sub G$ a compact, open subgroup
such that $\alpha(U)\sub U$.
The {\em calibration of~$G$ associated
with $\alpha$ and $U$} is the function
$\kappa\!: G\to \Z\cup\{\infty\}$ defined via
\begin{equation}\label{eqncalibration}
\kappa(x)\; :=\,
\left\{
\begin{array}{cl}
\infty &\;\mbox{if $\,x=1$;}\\
n & \;\mbox{if $\,x\in \alpha^n(U)\setminus \alpha^{n-1}(U)$.}
\end{array}
\right.
\end{equation}
\end{defn}
\begin{numba}\label{theset}
Let $G$ be a $C^1$-Lie group over $\K=\F_q(\!(X)\!)$
and $\alpha\!: G\to G$ be a contractive $C^1$-automorphism of~$G$
such that $s_G(\alpha^{-1})=q$. Let $\phi\!: P \to Q \sub L(G)$
be a chart for $G$ around~$1$ such that $\phi(1)=0$ and
$T_1(\phi)= \id_{L(G)}$. Let $\|.\|$ be an ultrametric norm on $L(G)$.
For each $r>0$, let $B_r$ be the ball of radius~$r$ in $(L(G),\|.\|)$;
let $\rho>0 $ such that $B_\rho \sub Q$, and set
$W_r:=\phi^{-1}(B_r)$ for $0<r\leq \rho$.
Let $U$ be a compact open subgroup of $G$
such that $\alpha(U)\sub U$.
\end{numba}
\begin{la}\label{prepcalib}
In the situation of {\bf\ref{theset}},
there exist $N\in \N_0$
and real numbers $a\leq b$ such that $q^b\leq \rho$ and
\begin{equation}\label{calib1}
W_{q^{a-n/d}}\; \sub \; \alpha^{n+N}(U)\; \sub W_{q^{b-n/d}}\qquad\mbox{for all $\,n\in \N_0$.}
\end{equation}
\end{la}
\begin{proof}
Let $\lambda_1,\ldots, \lambda_d$ (where $d:=\dim_\K \! L(G)$)
be the eigenvalues of $L(\alpha^{-1})$ in $\wb{\K}$,
each of which has absolute value $>1$ as a consequence of Lemma~\ref{wilcontra}\,(d).
Since $G$ has small subgroups tidy for $\alpha$ (Lemma~\ref{wilcontra}\,(c)), it also has small subgroups
tidy for~$\alpha^{-1}$. Hence
\begin{equation}\label{scaalpinv}
q\, =\, s_G(\alpha^{-1})\;\,  =\; \prod_{\stackrel{j\in\{1,\ldots, d\}}{{\scriptscriptstyle \text{s.t.} \, |\lambda_j|\geq 1}}} |\lambda_j|
\;\,=\;
\prod_{j=1}^d \: |\lambda_j|\,,
\end{equation}
using the main result of~\cite{POS} to obtain the second equality.
Now Lemma~\ref{crucialla} shows that
$|\lambda_j|=q^{1/d}$ for each $\,j\in \{1,\ldots, d\}$,
and it provides an ultrametric norm $\|.\|'$
on $L(G)$ such that
$\|L(\alpha^{-1}).x\|=q^{1/d}\|x\|'$
for each $x\in L(G)$. Hence
\[
\|\beta(x)\|'\; =\; q^{-1/d}\, \|x\|'\, , \qquad\mbox{for each $\,x\in L(G)$,}
\]
where $\beta:=L(\alpha)$.
Set $B_r':=\{x\in L(G)\!: \|x\|'<r\}$ for $r>0$.
There exists an open $0$-neighbourhood $R\sub Q$
such that $\alpha(\phi^{-1}(R))\sub P$. Then
\[
\gamma\!: R\to Q\, ,\quad \gamma(x)\, :=\, \phi(\alpha(\phi^{-1}(x)))
\]
expresses $\alpha$ in local coordinates.
We have $\gamma'(0)=T_1(\phi) L(\alpha) T_1(\phi)^{-1}=L(\alpha)=\beta$
as $T_1(\phi)= \id_{L(G)}$. Thus $\gamma'(0)=\beta$ is invertible.
Now the Inverse Function Theorem
in the quantitative form \cite[Prop.\,7.1\,(b)$'$]{IMP}
shows that there exists $\delta>0$ with $B'_\delta\sub R$
such that
$\gamma(B_r')=\beta(B_r')=B_{rq^{-1/d}}'$ for each $r\in \;]0,\delta]$
and hence
\begin{equation}\label{primballs}
\gamma^n(B'_r)=B_{rq^{-n/d}}'\qquad\mbox{for all $\,r\in \; ]0,\delta]$ and $\,n\in \N_0$.}
\end{equation} 
All ultrametric norms on the finite-dimensional $\K$-vector space $L(G)$ being equivalent \cite[Thm.\,13.3]{Sch},
there exists $c\geq 1$ such that $c^{-1}\|.\|\leq \|.\|'\leq c \|.\|$
and thus
$B_{c^{-1}r}\sub B_r'\sub B_{cr}$
for each $r>0$.
Set $v:=\min\{1,\, \frac{\rho}{c \delta}\}$.
Then $B'_{v\delta}\sub B'_\delta\sub R$. Choose $N\in \N$ such that
$\alpha^N(U)\sub \phi^{-1} (B'_{v\delta})$.
There exists $u\in \;]0,v]$ such that
$B'_{u\delta}\sub \phi(\alpha^N(U))$.
Then
\[
B'_{u\delta q^{-n/d}}\;\sub\;
\gamma^n(\phi(\alpha^N(U)))\;\sub \; B'_{v\delta q^{-n/d}}\qquad\mbox{for each $\, n\in \N_0$,}
\]
by~(\ref{primballs}).
Therefore $B_{c^{-1}u\delta q^{-n/d}}\sub \gamma^n(\phi(\alpha^N(U)))\sub B_{cv\delta q^{-n/d}}$.
Using that\linebreak
$\gamma^n(\phi(\alpha^N(U)))=\phi(\alpha^{n+N}(U))$, we obtain
\begin{equation}\label{localform}
B_{c^{-1}u\delta q^{-n/d}}\; \sub \; \phi(\alpha^{n+N}(U))\; \sub \; B_{cv\delta q^{-n/d}}
\qquad\mbox{for each $\,n\in \N_0$.}
\end{equation}
Note that $cv\delta q^{-n/d}\leq c v \delta\leq \frac{c \delta \rho}{c\delta}=\rho$
be definition of~$v$. Applying $\phi^{-1}$ to (\ref{localform}),
we see that (\ref{calib1}) holds with $a:=\log_q\big(\frac{u\delta}{c}\big)$ and $b:=\log_q (cv\delta)$.
\end{proof}
In the following lemma, we retain the situation of Lemma~\ref{prepcalib};
in particular, $d:=\dim_\K L(G)$.
For $\|.\|$ as in {\bf\ref{theset}}, define
$\mu\!: L(G)\to\R\cup\{\infty\}$,
$\mu(x):=-\log_q \|x\|$ (where $\log_q(0):=-\infty$).
Thus $\|x\|=q^{-\mu(x)}$.
The calibration
$\kappa \!:  G \to \Z \cup\{\infty\}$ associated with $\alpha$ and $U$
(as in (\ref{eqncalibration})) can be used to approximate~$\mu$ on a $0$-neighbourhood,
up to a bounded error:
\begin{la}[Calibration Lemma]\label{canibal}
There exists $R\in [0,\infty[$ such that
\begin{equation}\label{calib2}
\mu(\phi(x))d-R \;\leq \; \kappa(x)\;\leq \; \mu(\phi(x)) d + R\, ,
\quad
\mbox{for each $\,x\in \alpha^N(U)$.}
\end{equation}
%
\end{la}
\begin{proof}
Let $x\in \alpha^N(U)\setminus\{1\}$. Then $N\leq \kappa(x)<\infty$,
entailing that $n:=\kappa(x)-N\in \N_0$.
Using (\ref{calib1}), we obtain
\[
\begin{array}{ccccclll}
x & \in & \alpha^{\kappa(x)}(U) & = & \alpha^{n+N}(U) & \sub & W_{q^{b-n/d}} & \mbox{and}\\
x &\not\in & \alpha^{\kappa(x)+1}(U) &=& \alpha^{n+1+N}(U) &\supseteq & W_{q^{a-(n+1)/d}}\,. &
\end{array}
\]
Hence $q^{a-(n+1)/d}\leq \|\phi(x)\|< q^{b-n/d}$ and thus $(n+1)/d-a\geq \mu(\phi(x))> n/d-b$,
whence
$n+ 1-a d \geq \mu(\phi(x)) d>n-bd$
and therefore
\[
\kappa(x)-N+1-ad \;\geq\; \mu(\phi(x)) d\; >\; \kappa(x)-N-bd\,.
\]
Hence the assertion holds with $R:=
\max\{0,\, 1-ad-N,\,  N+bd\}$.
\end{proof}
The following
immediate consequence of Lemma~\ref{canibal}
is what we actually need:
\begin{prop}\label{corcalib}
Let $G$ and $G^*$ be
$C^1$-Lie groups over~$\K$,
of dimensions $d$ and $d^*$, resp.,
such that $G=G^*$
as a topological group.
Let $\phi\!: P\to Q\sub L(G)$
and $\phi^*\!: P^*\to Q^*$
be charts of $G$ and $G^*$ around~$1$
taking $1$ to~$0$ and such that $T_1(\phi)=\id_{L(G)}$
and $T_1(\phi^*)=\id_{L(G^*)}$.
Let $\|.\|$ and $\|.\|_*$ be ultrametric norms
on $L(G)$ and $L(G^*)$, respectively.
If there exists a contractive automorphism
$\alpha$ of $G$
which is a $C^1$-automorphism for both $G$ and $G^*$
and such that $s_G(\alpha^{-1})=q$,
then there exists an identity neighbourhood $W\sub P\cap P^*$
and positive real numbers $a<b$ such that
\[
a\, \|\phi(x)\|^{d/d^*} \;\leq \; \|\phi^*(x)\|_*\; \leq \; b\, \|\phi(x)\|^{d/d^*}\qquad \mbox{for each $\,x\in W$.}
\]
\end{prop}
$\;$\\[-12.3mm]
$\;$\vspace{6mm}\Punkt

\noindent
It may very well happen that $d\not=d^*$ in the situation of
Proposition~\ref{corcalib}, and in fact any combination
of positive integers $(d,d^*)\in \N^2$ can occur, as the following
example shows.
\begin{example}\label{defntau}
Let $G:=(\K,+)$, equipped with
its usual one-dimensional $\K$-analytic Lie group structure.
Then the (linear) map
\[
\tau\!: \K \to \K \, , \qquad \tau(z)\, :=\, X \cdot z
\]
apparently is a $\K$-analytic automorphism of~$G$.
The isomorphism of topological groups
\[
\phi^*:\, \K\to \K^2\, ,\qquad \sum_j a_jX^j\mto
\left( \sum_j a_{2j}X^j,\, \sum_j a_{2j+1}X^j\right)
\]
can be used as a global chart for $\K$ turning it into
a $2$-dimensional $\K$-analytic Lie group~$G^*$
isomorphic to $(\K^2,+)$.
Note that $\tau$ also is a $\K$-analytic automorphism of~$G^*$,
because $\phi^*\circ \tau\circ (\phi^*)^{-1}\!: \K^2\to\K^2$
is the $\K$-linear map
\[
\left(\begin{array}{c}
u\\
v
\end{array}
\right)
\; \mto\;
\left(
\begin{array}{cc}
0 & X\\
1 & 0
\end{array}
\right)\,
\left(\begin{array}{c}
u\\
v
\end{array}
\right)\,.
\]
By the following lemma, $\tau$ is contractive and $s_G(\tau^{-1})=q$.
Hence, we are in the situation of Proposition~\ref{corcalib}.\\[3mm]
{\em Variants}: Given any positive integer $d^*$, we obtain an
isomorphism of topological groups $\K\isom \K^{d^*}$
if we replace $2\Z$ and $2\Z+1$ with $d^*\Z$ and its cosets
in the preceding construction. This enables
us to turn $\K$ into a $\K$-analytic Lie group $G^*$ of dimension~$d^*$,
such that $\tau$ is a $\K$-analytic automorphism of~$G^*$.
\end{example}
We used the following simple fact:
\begin{la}\label{tau}
$\tau$ is a contractive automorphism of $(\K,+)$ with $s_\K(\tau^{-1})=q$.
\end{la}
\begin{proof}
Because $|\tau^n(z)|=|X^nz|=q^{-n}|z|$ for each $z\in\K$ and $n\in \N$,
clearly $\tau$ is contractive.
Since $\tau(\bO)=X\bO\sub \bO$,
Lemma~\ref{wilcontra}\,(b)
shows that
$\bO$ is tidy for~$\tau$,
with $\bO_+=\{0\}$ and $\bO_- =\bO$.
Therefore $s_\K(\tau^{-1})=[\tau^{-1}(\bO)\, :\, \bO]=[X^{-1}\bO\,:\, \bO]=|\F_q|=q$.
\end{proof}
\section{Examples of non-analytic {\boldmath $C^\infty$\/}-Lie groups}\label{secnanagp}
In this section, we present examples of smooth
$\K$-Lie groups which cannot be turned into $\K$-analytic Lie groups.
\begin{numba}\label{defgreat}
Let $\lambda=(\ell_j)_{j\in \N_0}$,
$\alpha$ and $\beta$
be as in {\bf\ref{defalpbet}}.
Every element $z=\sum_{j=-n}^\infty a_jX^j\in \K$ (where $n\in \N$)
can be written
as $z=z'+z''$ with $z':=\sum_{j=-n}^{-1}a_jX^j$
and $z'':=\sum_{j=0}^\infty a_jX^j$.
We define 
\[
\bar{\beta}\!: \K \to \K \, , \qquad \bar{\beta}(z)\, :=\, \beta(z'')\quad\qquad\mbox{and}
\]
\[
\!\!\!\!\!\!\!\!\!\!\!\!\!\bar{\alpha}\!: \K \to \K \,  , \qquad \bar{\alpha}(z)\, := z'+ \alpha(z'')\,.
\]
Thus $\bar{\alpha}|_\bO=\alpha$, $\bar{\beta}|_\bO=\beta$,
and $\bar{\alpha}=\id_\K +\bar{\beta}$. Clearly $\bar{\beta}$
is a continuous homomorphism
and $\bar{\alpha}$ a bicontinuous automorphism of~$(\K,+)$.
We define
\[
\tau\!: \K \to \K \, , \qquad \tau(z)\, :=\, X \cdot z
\]
and recall from {\bf \ref{defntau}} and Lemma~\ref{tau}
that~$\tau$ is a $\K$-analytic, contractive automorphism of $(\K,+)$
such that $s_\K(\tau^{-1})=q$.
We equip
\begin{equation}\label{thegroup}
G \; :=\; \K\dsemi \langle \bar{\alpha},\tau\rangle
\end{equation}
with the uniquely determined topology making it a topological group
with $\K$ (in its usual topology) as an open, normal subgroup;
here $\langle \bar{\alpha},\tau\rangle$ denotes the subgroup
of $\Aut(\K,+)$
generated by $\bar{\alpha}$ and $\tau$.
\end{numba}
\begin{thm}\label{mainthm}
For $\lambda=(\ell_j)_{j\in \N_0}$ as in Proposition~{\rm\ref{mainwork2}},
the topological group $G$ described in {\rm (\ref{thegroup})}
has the following properties:
\begin{itemize}
\item[\rm (a)]
$G$ admits a $1$-dimensional smooth $\K$-Lie group structure
compatible with the given topological group structure.
\item[\rm (b)]
$G$ does not admit a $\K$-analytic Lie group structure
compatible with the given topological group structure.
\end{itemize}
\end{thm}
{\bf Proof.} (a)
We equip $(\K,+)$ with its usual $1$-dimensional
$C^\infty$-Lie group structure.
Since the restriction of $\bar{\alpha}$
to the open subset $\bO\sub \K$ coincides with $\alpha$,
which is a $C^\infty$-automorphism of~$\bO$
by Proposition~\ref{mainwork2},
we see that $\bar{\alpha}$ is a $C^\infty$-automorphism of $(\K,+)$.
Apparently also $\tau$ is a $C^\infty$-automorphism of $(\K,+)$,
and hence so is each $\gamma\in \langle \bar{\alpha},\tau\rangle$.
Now
standard arguments provide
a (uniquely determined) $C^\infty$-manifold structure on~$G$
making it a smooth Lie group and
such that $\K$, with its given smooth manifold structure,
is an open submanifold of~$G$ (see \cite[Prop.\,1.18]{ZOO}).\vspace{2mm}

(b) The proof has two parts,
the first of which can be re-used later.
\begin{numba}\label{reusable}
Let us assume, to the contrary, that
there exists a $\K$-analytic manifold structure on the group~$G$
making it a $\K$-analytic Lie group~$G^*$
and compatible with the given topology.
We let $L(G^*)$ be the Lie algebra of $G^*$
and set $d:=\dim_\K L(G^*)\in \N$.
The open subgroup $\K$ of $G^*$
inherits a $\K$-analytic Lie group structure from~$G^*$.
We write~$H$ for~$\K$, equipped with this $d$-dimensional
$\K$-analytic Lie group
structure. Thus $L(H)=L(G^*)$,
and both $\tau$ and $\bar{\alpha}$ are $\K$-analytic automorphisms
of~$H$. We let $\phi\!: P\to Q\sub L(H)$ be a chart for~$H$ around~$0$
such that $\phi(0)=0$ and
$T_0(\phi)=\id_{L(H)}$.
Furthermore, we choose an ultrametric norm
$\|.\|$ on $L(H)$.
Using Proposition~\ref{corcalib},
we find $N\in \N_0$ and positive real numbers $a<b$ such that $X^N\bO\sub P$
and
\begin{equation}\label{estm1}
a\, \|\phi(z)\|^d\;\leq\; |z|\;\leq\; b\, \|\phi(z)\|^d\qquad\mbox{for each \,$z\in X^N\bO$.}
\end{equation}
Set $W:=\phi(X^N\bO)$.
Since $\beta(X^N\bO)\sub X^N\bO$, we can define
$\gamma\!: W\to W$, $\gamma(w):=\phi(\beta(\phi^{-1}(w)))$.
\end{numba}
\begin{numba}
Using Lemma~\ref{givnonana}, we want to show that $\gamma$ is not analytic.
Given $n\in \N$, there exists $k\geq N$ such that $\ell_j\geq (n+1)j$
for each $j\geq k$ and hence
\begin{equation}\label{estm2}
|\beta(z)|\; \leq \; |z|^{n+1}\qquad\mbox{for each $\,z\in X^k\bO$.}
\end{equation}
For $w\in W$ such that $\|w\|\leq q^{-k/d}b^{-1/d}$,
we have $|\phi^{-1}(w)|\leq b\|w\|^d\leq q^{-k}$ (using (\ref{estm1}))
and hence
$|\beta(\phi^{-1}(w))|\leq |\phi^{-1}(w)|^{n+1}$, by~(\ref{estm2}).
As a consequence,
\[
\|\gamma(w)\|\,=\, \|\phi(\beta(\phi^{-1}(w)))\|\,\leq\,
a^{-1/d}|\phi^{-1}(w)|^{(n+1)/d}\,\leq \, a^{-1/d}b^{(n+1)/d}\|w\|^{n+1}\,,
\]
entailing that $\|\gamma(w)\|/\|w\|^n\to 0$ as $w\to 0$.
Since~$\gamma$ is injective,
Lemma~\ref{givnonana}
shows that $\gamma$
cannot be analytic.
As a consequence, $\bar{\beta}\!: H\to H$ is not analytic either,
nor is $\bar{\alpha}=\id_H+\bar{\beta}\!: H\to H$.
But $\bar{\alpha}\!: H\to H$ has to be analytic,
being the restriction of an inner automorphism of the analytic Lie group
$G^*$ to its open, normal subgroup~$H$.
We have reached a contradiction.\Punkt
\end{numba}
\begin{rem}\label{linus}
In particular, Theorem~\ref{mainthm} shows that
there exist smooth Lie groups over local fields of positive characteristic
without a {\em smoothly\/} compatible analytic Lie group
structure.
\end{rem}
\section{\!\!\!Examples of {\boldmath $C^n$}-Lie groups that are not {\boldmath $C^{n+1}$}}\label{secnCngp}
A variant of the proof of Theorem~\ref{mainthm} shows the following:
\begin{thm}\label{mainthm2}
Given $n\in \N$, let $\theta
$
and $\lambda=(\ell_j)_{j\in \N_0}$
be as in Proposition~{\rm\ref{mainwork1}}.
Then the topological group~$G$ from {\rm (\ref{thegroup})}
has the following properties:
\begin{itemize}
\item[\rm (a)]
$G$ admits a $1$-dimensional $\K$-Lie group structure
of class~$C^n$
compatible with the given topological group structure.
\item[\rm (b)]
$G$ does not admit a $\K$-Lie group structure
of class $C^{n+1}$
compatible with the given topological group structure.
\end{itemize}
\end{thm}
\begin{proof}
The proof of Theorem~\ref{mainthm}\,(a)
and the first half of the proof of Theorem~\ref{mainthm}\,(b)
(viz.\ {\bf\ref{reusable}})
carry over to the present situation, if we replace the word ``smooth''
with ``$C^n$'' and ``$\K$-analytic'' with~``$C^{n+1}$'' there.
Thus Part\,(a) of the present theorem holds. To complete the proof of (b),
we retain the notations from
{\bf\ref{reusable}}. Using Lemma~\ref{nonCn}, we want to show
that~$\gamma$ is not~$C^{n+1}$ (whence neither is $\bar{\beta}$,
nor $\bar{\alpha}$, which is absurd).
By Lemma~\ref{lamgivesas}, there exist positive real numbers $A<B$
such that
\begin{equation}\label{finaleq}
A|z|^\theta\;\leq \; |\beta(z)|\;\leq\; B|z|^\theta,\qquad \mbox{for each $\, z\in \bO$.}
\end{equation}
For each $w\in W$, we then have
\begin{eqnarray*}
\|\gamma(w)\| & = & \|\phi(\beta(\phi^{-1}(w)))\| \;\leq \;
a^{-1/d}|\beta(\phi^{-1}(w))|^{1/d}\\[1mm]
&\leq & a^{-1/d}B^{1/d}|\phi^{-1}(w)|^{\theta/d}\; \leq \;
a^{-1/d}b^{\theta/d}B^{1/d}\|w\|^\theta\qquad\mbox{and}
\end{eqnarray*}
\begin{eqnarray*}
\!\! \|\gamma(w)\| &=& \|\phi(\beta(\phi^{-1}(w)))\|\;\geq \;
b^{-1/d}|\beta(\phi^{-1}(w))|^{1/d}\qquad\qquad\quad\\[1mm]
&\geq & b^{-1/d}A^{1/d}|\phi^{-1}(w)|^{\theta/d}\;
\geq\; b^{-1/d}A^{1/d}a^{\theta/d}\|w\|^\theta\,,
\end{eqnarray*}
using (\ref{estm1}) and (\ref{finaleq}).
Set $C:=b^{-1/d}A^{1/d}a^{\theta/d}$
and $D:=a^{-1/d}b^{\theta/d}B^{1/d}$.
Then
\[
C\|w\|^\theta\;\leq\; \|\gamma(w)\|\;\leq\; D\|w\|^\theta\qquad
\mbox{for each $\,w\in W$,}
\]
whence Lemma~\ref{nonCn} shows that $\gamma$ is not $C^{n+1}$, contradiction.
\end{proof}
The following problems remain open:
\begin{problem}
{\rm Does every smooth $\K$-Lie group
have an open subgroup which admits a smoothly compatible
analytic Lie group structure\,?}
\end{problem}
\begin{problem}
{\rm Given $n\in \N$,
does every $\K$-Lie group of class $C^n$
have an open subgroup which admits a $C^n$-compatible
$\K$-Lie group structure of class~$C^{n+1}$\,?}
\end{problem}
%
%
%
{\em Acknowledgements.}
The present research
was carried out at the University of Newcastle (N.S.W.)
in August 2004,
supported by DFG grant 447 AUS-113/22/0-1
and ARC grant LX 0349209.
The author profited from
discussions with
George Willis 
and Udo Baumgartner.
Remark~\ref{linus} answers a question
posed to the author by Linus Kramer (Darmstadt) in 2003.
{\footnotesize

Helge Gl\"{o}ckner, TU Darmstadt, Fachbereich Mathematik AG~5,
Schlossgartenstr.\,7,\\
64289 Darmstadt, Germany. E-Mail: gloeckner\at{}mathematik.tu-darmstadt.de}
\end{document}